\documentclass[12pt]{article}
\usepackage{amsmath}
\usepackage{amssymb}
\usepackage{amsthm}
\usepackage{mathrsfs}
\usepackage[colorlinks,linkcolor=blue,citecolor=blue]{hyperref}

\allowdisplaybreaks[0]

\makeatletter \@addtoreset{equation}{section}

\makeatletter\@addtoreset{figure}{section}\makeatother

\newtheorem{thm}{Theorem}[section]
\newtheorem{lem}[thm]{Lemma}
\newtheorem{prop}[thm]{Proposition}

\begin{document}
\begin{center}
{\large \bf Sun's log-concavity conjecture\\
 on the Catalan-Larcombe-French sequence}
\end{center}
\begin{center}
James J.Y. Zhao\\[6pt]

Center for Applied Mathematics\\
Tianjin University, Tianjin 300072, P. R. China\\[8pt]

Email: {\tt jjyzhao@tju.edu.cn}
\end{center}

\noindent\textbf{Abstract.}
Let $\{P_n\}_{n\geq 0}$ denote the Catalan-Larcombe-French sequence, which naturally came up  from the series expansion of the complete elliptic integral of the first kind. In this paper, we prove
the strict log-concavity of the sequence $\{\sqrt[n]{P_n}\}_{n\geq 1}$, which was originally conjectured by Sun. We also obtain the strict log-concavity of the sequence $\{\sqrt[n]{V_n}\}_{n\geq 1}$, where $\{V_n\}_{n\geq 0}$ is the Fennessey-Larcombe-French sequence arising in the series expansion of the complete elliptic integral of the second kind.

\noindent \emph{AMS Classification 2010:} Primary 05A20


\section{Introduction}

The main objective of this paper is to prove the log-concavity of the sequence $\{\sqrt[n]{P_n}\}_{n\geq 1}$ and the sequence $\{\sqrt[n]{V_n}\}_{n\geq 1}$, where $\{P_n\}_{n\geq 0}$ and $\{V_n\}_{n\geq 0}$, known as the Catalan-Larcombe-French sequence and the Fennessey-Larcombe-French sequence respectively, are given by
\begin{align}
(n+1)^2P_{n+1}&=8(3n^2+3n+1)P_n-128n^2P_{n-1},\label{eq-rclf}\\
n(n+1)^2 V_{n+1}&=8n(3n^2+5n+1)V_n-128(n-1)(n+1)^2V_{n-1},\label{eq-rec}
\end{align}
with the initial values $P_0=V_0=1$ and $P_1=V_1=8$. These two sequences came up naturally from the series expansions of the complete elliptic integrals. For more information on $\{P_n\}_{n\geq 0}$ and
$\{V_n\}_{n\geq 0}$, see \cite{Catalan, jv2010, lf2000, lff2002, lff2003}.

Let us first give an overview of some background. Recall that a sequence $\{a_n\}_{n \geq 0}$ of real numbers is said to be log-concave (resp. log-convex) if
$$a_{n}^2\geq a_{n-1}a_{n+1} \mbox{  (resp. $a_{n}^2\leq a_{n-1}a_{n+1}$)}$$
for all $n\geq 1$, and it is strictly log-concave (resp. strictly log-convex) if the inequality is strict.
Clearly, for some integer $N>0$, the positive sequence $\{a_n\}_{n\geq N}$ is log-concave (resp. log-convex) if and only if the sequence $\{a_{n+1}/a_n\}_{n\geq N}$ is decreasing (resp. increasing).

The strict log-concavity of $\{\sqrt[n]{P_n}\}_{n\geq 1}$ was first conjectured by Sun \cite{sun2013}, who also conjectured that the sequences $\{P_{n+1}/P_{n}\}_{n\geq 0}$ and $\{\sqrt[n]{P_n}\}_{n\geq 1}$ are strictly increasing. In fact, motivated by Firoozbakht's conjecture on the strictly decreasing property of the sequence $\{\sqrt[n]{p_n}\}_{n\geq 1}$, where $p_n$ is the $n$-th prime number, see \cite[p.\!\! 185]{Ribenboim}, Sun \cite{sun2013} studied the monotonicity of many number theoretical and combinatorial sequences.

The strictly increasing property of the sequence $\{P_{n+1}/P_{n}\}_{n\geq 0}$ and $\{\sqrt[n]{P_n}\}_{n\geq 1}$ has been confirmed by Xia and Yao \cite{XY2013}, and independently by Zhao \cite{zhao2014}. We would like to point out that Zhao \cite{zhao2014} only proved the
log-convexity of $\{{P_n}\}_{n\geq 0}$, which also implies the monotonicity of $\{\sqrt[n]{P_n}\}_{n\geq 1}$ by using a result due to Wang and Zhu \cite{WangZhu}.

For a positive sequence $\{a_n\}_{n\geq 0}$ satisfying a three-term recurrence relation, Chen, Guo and Wang \cite{CGW2014} obtained a useful criterion to determine
the log-concavity of $\{\sqrt[n]{a_n}\}_{n\geq N}$ for some positive integer $N$. While, their criterion does not apply to the
the Catalan-Larcombe-French sequence, and Sun's conjecture on the log-concavity of $\{\sqrt[n]{P_n}\}_{n\geq 1}$ remains open. The first main result of this paper is as follows.

\begin{thm}\label{th-rlcc}
The sequence $\{\sqrt[n]{P_n}\}_{n\geq 1}$ is strictly log-concave, that is, for $n\geq 2$,
\begin{align}\label{eq-rlcc}
\left(\sqrt[n]{P_n}\right)^2>\sqrt[n-1]{P_{n-1}}\cdot\sqrt[n+1]{P_{n+1}}.
\end{align}
\end{thm}

Since the Fennessey-Larcombe-French sequence $\{V_n\}_{n\geq 0}$ is closely related to
the sequence $\{P_n\}_{n\geq 0}$, we are led to study the
log-behavior of $\{\sqrt[n]{V_n}\}_{n\geq 1}$. The second main result of this paper is as follows.
\begin{thm}\label{th-nrlc}
The sequence $\{\sqrt[n]{V_n}\}_{n\geq 1}$ is strictly log-concave, that is, for $n\geq 2$,
\begin{align}\label{eq-nrlc}
\left(\sqrt[n]{V_n}\right)^2>\sqrt[n-1]{V_{n-1}}\cdot\sqrt[n+1]{V_{n+1}}.
\end{align}
\end{thm}

This paper is organized as follows. In Section \ref{s-P}, we give a proof of Theorem \ref{th-rlcc} by establishing a lower bound and an upper bound for the ratio $P_n/P_{n-1}$. The proof of Theorem \ref{th-nrlc} is similar to that of  Theorem \ref{th-rlcc}, which will be given in Section \ref{s-V}. Finally, we give a new proof of the monotonicity of $\{\sqrt[n]{P_n}\}_{n\geq 1}$ by using Theorem \ref{th-rlcc}, and obtain the monotonicity of $\{\sqrt[n]{V_n}\}_{n\geq 1}$ by using Theorem \ref{th-nrlc}.

\section{Proof of Theorem \ref{th-rlcc}}\label{s-P}

In this section we aim to prove Theorem \ref{th-rlcc}.
To this end, we first establish a lower bound and an upper bound for the ratio $P_n/P_{n-1}$ that will lead to the strict log-concavity of the sequence $\{\sqrt[n]{P_n}\}_{n\geq 1}$.

\begin{lem}\label{lem-p}
For any integer $n\geq 1$, let
\begin{align}\label{eq-fn}
f(n)=\frac{16(n-1)}{n}.
\end{align}
Then for $n\geq 5$, we have
\begin{align}\label{eq-bP}
f(n-1)<\frac{P_n}{P_{n-1}}<f(n).
\end{align}
\end{lem}

\proof
For notational convenience, let $v(n)=P_n/P_{n-1}$. We first use induction on $n$ to prove
$v(n)>f(n-1)$ for $n\geq 5$. By the recurrence \eqref{eq-rclf}, we have
\begin{align}\label{eq-Prr}
v(n+1)=\frac{8(3n^2+3n+1)}{(n+1)^2}-\frac{128n^2}{(n+1)^2 v(n)},\quad n\geq 1,
\end{align}
with the initial value $v(1)=8$.
Clearly, $v(5)=2152/169>12=f(4)$. Assume that $v(n)>f(n-1)$, and we proceed to prove that $v(n+1)>f(n)$ for $n\geq 5$. Observe that
\begin{align*}
 v(n+1)-f(n)
=&\ \frac{8(3n^2+3n+1)}{(n+1)^2}-\frac{128n^2}{(n+1)^2 v(n)}-\frac{16(n-1)}{n}\\[5pt]
=&\ \frac{8(n^3+n^2+3n+2)v(n)-128n^3}{n(n+1)^2 v(n)}.
\end{align*}
By the induction hypothesis, we have $v(n)>f(n-1)>0$ for $n\geq 5$, therefore
\[  v(n+1)-f(n)
>\frac{8(n^3+n^2+3n+2)f(n-1)-128n^3}{n(n+1)^2 v(n)}
=\frac{128(n^2-4n-4)}{(n-1)n(n+1)^2 v(n)}
>0
\]
for $n\geq 5$, since $n^2-4n-4=(n+1)(n-5)+1>0$ for $n\geq 5$.
This proves that $v(n)>f(n-1)$ for $n\geq 5$.

The inequality $v(n)<f(n)$ for $n\geq 5$ can be obtained in the same manner,
and the detailed proof is omitted here.
This completes the proof.
\qed

With the bounds given in Lemma \ref{lem-p} we are now able to prove Theorem \ref{th-rlcc}.

\noindent{\it Proof of Theorem \ref{th-rlcc}.}
Note that for $n\geq 2$, the inequality \eqref{eq-rlcc} can be rewritten as
\begin{align*}
\frac{\sqrt[n]{P_{n}}}{\sqrt[n-1]{P_{n-1}}}
>\frac{\sqrt[n+1]{P_{n+1}}}{\sqrt[n]{P_n}},
\end{align*}
or equivalently,
\begin{align}\label{eq-rpn}
\left(\frac{P_n}{P_{n-1}}\right)^{n(n+1)}
>P_n^2\left(\frac{P_{n+1}}{P_n}\right)^{n(n-1)}.
\end{align}

By the recurrence \eqref{eq-rclf}, it is easy to verify that \eqref{eq-rpn} holds for $2\leq n\leq 6$. We proceed to prove that \eqref{eq-rpn} is true for $n\geq 7$. By Lemma \ref{lem-p} we have
\begin{align*}
\left(\frac{P_n}{P_{n-1}}\right)^{n(n+1)}
 \hskip -0.5ex -P_n^2\left(\frac{P_{n+1}}{P_n}\right)^{n(n-1)}\hskip -0.5ex
>& (f(n-1))^{n(n+1)}-P_n^2(f(n+1))^{n(n-1)}\\[3pt]
=& \left(\frac{16(n-2)}{n-1}\right)^{n(n+1)}
    -P_n^2\left(\frac{16n}{n+1}\right)^{n(n-1)}\\[3pt]
=& 16^{n(n-1)}\left(16^{2n}\left(\frac{n-2}{n-1}\right)^{n(n+1)}
    -P_n^2\left(\frac{n}{n+1}\right)^{n(n-1)}\right).
\end{align*}
It suffices to prove that
\begin{align*}
16^{2n}\left(\frac{n-2}{n-1}\right)^{n(n+1)}
-P_n^2\left(\frac{n}{n+1}\right)^{n(n-1)}>0,
\end{align*}
or equivalently
\begin{align*}
\frac{P_n^2}{16^{2n}}<\left(\frac{(n-2)(n+1)}{(n-1)n}\right)^{n(n-1)}
  \left(\frac{n-2}{n-1}\right)^{2n}.
\end{align*}
Thus we only need to show that, for $n\geq 7$,
\begin{align}\label{eq-sep}
\frac{P_n}{16^n}<\left(\frac{(n-2)(n+1)}{(n-1)n}\right)^{\frac{n(n-1)}{2}}
  \left(\frac{n-2}{n-1}\right)^{n}.
\end{align}
Let $l_n$ denote the term on the left hand side, and $r_n$ denote the term on the right hand side.
We claim that
\begin{itemize}
\item[(i)] the sequence $\{l_n\}_{n\geq 5}$ is strictly decreasing, and
\item[(ii)] the sequence $\{r_n\}_{n\geq 5}$ is strictly increasing.
\end{itemize}

By Lemma \ref{lem-p}, we see that
\[ 0<\frac{l_n}{l_{n-1}}=\frac{P_n/P_{n-1}}{16}<1\]
for $n\geq 5$, which implies (i).

We proceed to prove (ii). Note that
$$
r_n=\left(1-\frac{1}{{n\choose 2}}\right)^{n\choose 2}\cdot
 \left(1-\frac{1}{n-1}\right)^{n-1}\cdot \left(1-\frac{1}{n-1}\right).
$$
The increasing property of the sequence $\{r_n\}_{n\geq 5}$ immediately follows from the
well-known fact that the sequence
$\{\left(1-\frac{1}{n}\right)^{n}\}_{n\geq 1}$ is strictly increasing.

It is easy to verify that $l_7<r_7$. Combining (i) and (ii), we get that
$l_n<r_n$ for $n\geq 7$, namely \eqref{eq-sep} holds. This completes the proof.
\qed

\section{Proof of Theorem \ref{th-nrlc}}\label{s-V}
In this section, we complete the proof of Theorem \ref{th-nrlc} in a similar way of that of Theorem \ref{th-rlcc}. For this purpose, we need a lower bound and an upper bound for the ratio $V_n/V_{n-1}$. For integer $n\geq 2$, let
\begin{align}\label{eq-hn}
h(n)=\frac{16(n^3-n^2+1)}{n^3-n^2},
\end{align}
which was introduced by Yang and Zhao \cite{YangZhao} in their study of the log-concavity
of the sequence $\{V_n\}_{n\geq 1}$.
They \cite{YangZhao} showed that $h(n)$ is a lower bound for the ratio $V_n/V_{n-1}$, precisely,
\begin{align}\label{eq-Vnhn}
\frac{V_n}{V_{n-1}}> h(n),
\end{align}
for $n\geq 4$. We further show that $h(n-1)$ is an upper bound for the ratio $V_n/V_{n-1}$.

\begin{lem}\label{lem-2}
Let $h(n)$ be given by \eqref{eq-hn}. Then for $n\geq 11$, we have
\begin{align}\label{eq-Vnru1}
\frac{V_n}{V_{n-1}}<h(n-1).
\end{align}
\end{lem}

\noindent{\it Proof. }
Let $g(n)=V_n/V_{n-1}$, by \eqref{eq-rec} it is clear that
\begin{align}\label{eq-gr}
g(n+1)=\frac{8(3n^2+5n+1)}{(n+1)^2}-\frac{128(n-1)}{ng(n)},
\quad n\geq 1,
\end{align}
with $g(1)=8$.
The inequality \eqref{eq-Vnru1} can be proved inductively based on the recurrence \eqref{eq-gr}. The proof is  similar to that of Lemma \ref{lem-p}, and hence is omitted here.
\qed

We proceed to prove Theorem \ref{th-nrlc}.

\noindent{\it Proof of Theorem \ref{th-nrlc}.} Note that for $n\geq 2$, the inequality \eqref{eq-nrlc} can be rewritten as
\begin{align*}
\frac{\sqrt[n]{V_{n}}}{\sqrt[n-1]{V_{n-1}}}
>\frac{\sqrt[n+1]{V_{n+1}}}{\sqrt[n]{V_n}},
\end{align*}
or equivalently,
\begin{align}\label{eq-rgn}
\left(\frac{V_n}{V_{n-1}}\right)^{n(n+1)}
>V_n^2\left(\frac{V_{n+1}}{V_n}\right)^{n(n-1)}.
\end{align}

By the recurrence \eqref{eq-rec}, it is easy to verify that \eqref{eq-rgn} holds for $2\leq n\leq 9$. We proceed to prove that \eqref{eq-rgn} is true for $n\geq 10$.
By \eqref{eq-Vnhn} we have
\[ \frac{V_n}{V_{n-1}}>h(n),\] for $n\geq 4$. By Lemma \ref{lem-2}, we have
\[ \frac{V_{n+1}}{V_n}<h(n), \] for $n\geq 10$.
Then for $n\geq 10$, it follows that
\begin{align*}
\left(\frac{V_n}{V_{n-1}}\right)^{n(n+1)}
  -V_n^2\left(\frac{V_{n+1}}{V_n}\right)^{n(n-1)}
>&\ (h(n))^{n(n+1)}-V_n^2(h(n))^{n(n-1)}\\[3pt]
=&\ (h(n))^{n(n-1)}\left((h(n))^n+V_n\right)\left((h(n))^n-V_n\right).
\end{align*}
Clearly, both $(h(n))^{n(n-1)}$ and $\left((h(n))^n+V_n\right)$ are positive for $n\geq 10$.
Thus it suffices to show that for $n\geq 10$,
\begin{align}\label{eq-vuh}
(h(n))^n-V_n>0.
\end{align}

Note that the relation $(24)$ in \cite{lff2002} gave an upper bound of $V_n$, that is, for $n\geq 1$,
\begin{align}\label{eq-Vu1}
V_n<(2n+1)\binom{2n}{n}^2.
\end{align}
Sasv\'{a}ri \cite[Corollary 1]{Sasvari} showed that for $n\geq 1$,
\begin{align}\label{eq-ubs}
\binom{2n}{n}<\frac{4^n}{\sqrt{\pi n}}e^{-\frac{1}{8n}+\frac{1}{192n^3}}.
\end{align}
It is clear that for $n\geq 1$,
\[ 0<e^{-\frac{1}{8n}+\frac{1}{192n^3}}<1. \]
Then it follows from \eqref{eq-ubs} that for $n\geq 1$
\begin{align}\label{eq-nfu}
\binom{2n}{n}<\frac{4^n}{\sqrt{\pi n}}.
\end{align}
 Combining \eqref{eq-Vu1} and \eqref{eq-nfu}, for $n\geq 1$ we have
\begin{align}\label{eq-Vu2}
V_n<\frac{2n+1}{\pi n}16^n<16^n.
\end{align}
Note that for $n\geq 2$
\begin{align}\label{eq-hnn}
h(n)=\frac{16(n^3-n^2+1)}{n^3-n^2}>16.
\end{align}
Thus by \eqref{eq-Vu2} and \eqref{eq-hnn}, we obtain
\eqref{eq-vuh}.
This completes the proof.
\qed

\section{The monotonicity of $\{\sqrt[n]{P_n}\}_{n\geq 1}$ and $\{\sqrt[n]{V_n}\}_{n\geq 1}$}

In this section, we aim to derive the monotonicity of $\{\sqrt[n]{P_n}\}_{n\geq 1}$ from
the log-concavity of $\{\sqrt[n]{P_n}\}_{n\geq 1}$, and to derive the monotonicity of $\{\sqrt[n]{V_n}\}_{n\geq 1}$ from
the log-concavity of $\{\sqrt[n]{V_n}\}_{n\geq 1}$. The main result of this section is as follows.

\begin{prop}\label{th-nP2}
Both $\{\sqrt[n]{P_n}\}_{n\geq 1}$ and $\{\sqrt[n]{V_n}\}_{n\geq 1}$ are strictly increasing.
\end{prop}

\proof

Recall that for a real sequence $\{a_n\}_{n\geq 0}$ with positive numbers, it was shown that
\[
\liminf\limits_{n\rightarrow\infty}\frac{a_{n+1}}{a_n}
\leq \liminf\limits_{n\rightarrow\infty}\sqrt[n]{a_n},
\]
and
\[
\limsup\limits_{n\rightarrow\infty}\sqrt[n]{a_n}
\leq \limsup\limits_{n\rightarrow\infty}\frac{a_{n+1}}{a_n},
\]
see Rudin \cite[\S 3.37]{Rudin}.
These two inequalities imply a well-known criterion, that is, if $\lim_{n\rightarrow\infty}\frac{a_n}{a_{n-1}}=c$, then $\lim_{n\rightarrow\infty}\sqrt[n]{a_n}=c$, where $c$ is a real number.

Note that it was proved in \cite[Eq. (30)]{lff2002} that $$\lim_{n\rightarrow\infty}\frac{P_n}{P_{n-1}}=
\lim_{n\rightarrow\infty}\frac{V_n}{V_{n-1}}=16.$$
Then it follows that
\begin{align*}
\lim\limits_{n\rightarrow \infty} \sqrt[n]{P_{n}}=\lim\limits_{n\rightarrow \infty} \sqrt[n]{V_{n}}=16.
\end{align*}
Consequently,
\[
\lim\limits_{n\rightarrow \infty} {\sqrt[n+1]{P_{n+1}}}/{\sqrt[n]{P_n}}
=\lim\limits_{n\rightarrow \infty} {\sqrt[n+1]{V_{n+1}}}/{\sqrt[n]{V_n}}=1. \]

By Theorems \ref{th-rlcc} and \ref{th-nrlc}, we see that
both $\{\sqrt[n+1]{P_{n+1}}/\sqrt[n]{P_n}\}_{n\geq 1}$ and $\{\sqrt[n+1]{V_{n+1}}/\sqrt[n]{V_n}\}_{n\geq 1}$
are strictly decreasing.  Thus for all $n\geq 1$, we have
\[ {\sqrt[n+1]{P_{n+1}}}/{\sqrt[n]{P_n}}>1, \quad {\rm and}\quad {\sqrt[n+1]{V_{n+1}}}/{\sqrt[n]{V_n}}>1. \]
This completes the proof.
\qed

\noindent{\it Remark.}
By employing a criterion due to Wang and Zhu \cite[Theorem 2.1]{WangZhu}, Yang and Zhao \cite{YangZhao} showed that the sequence $\{\sqrt[n]{V_{n+1}}\}_{n\geq 1}$ is strictly decreasing.

\vskip 1mm
\noindent {\bf Acknowledgments.} This work was supported by the 863 Program and the National Science Foundation of China.

\end{document}